\documentclass{amsart}
\usepackage{amsmath}
\usepackage{amssymb}
\usepackage{amscd}
\usepackage{curves}
\usepackage{amsrefs}
\newtheorem{prob}{Problem}

\theoremstyle{definition}

\def \CPb {\overline{\mathbf{CP}}^{\,2}}
\def \CP {{\mathbf{CP}}^{\,2}} 
\def \R {\mathbf{R}}
\def \Z {\mathbf{Z}}
\def \Sig{\Sigma}

\def \PO {{\mathbf{CP}}^{\,1}}
\def \SS {S^2\times S^2}

\def \vp {\varphi}

\def \a {\alpha}
\def \b {\beta}
\def \g {\gamma}
\def \d {\delta}

\def \lam {\lambda}

\def \o {\omega}
\def \s {\sigma}

\def \bd {\partial}
\def \x {\times}
\def \- {\setminus}
\def \C {\subset}

\def \ch {\chi_{{}_h}}
\def \ssw {\text{SW}}
\def \sw {\mathcal{SW}}
\def \swb{\overline {\mathcal{SW}}}

\def \DD {\Delta}

\def\fs{\mathfrak{s}}

\def\hk{\widehat{k}}

\begin{document}

\baselineskip.5cm

\title{Will we ever classify simply-connected smooth 4-manifolds?}

\author[Ronald J. Stern]{Ronald J. Stern}
\address{Department of Mathematics, University of California, Irvine, CA 92697} \email{rstern@uci.edu}
\thanks{The author was partially supported by NSF Grant DMS0204041}


\begin{abstract}
These notes are adapted from two talks given at the 2004 Clay
Institute Summer School on {\it Floer homology, gauge theory, and low
dimensional topology} at the Alfred R\'enyi Institute. We will quickly
review what we do and do not know about the existence and uniqueness
of smooth and symplectic structures on closed, simply-connected
$4$-manifolds. We will then list the techniques used to date and
capture the key features common to all these techniques. We finish
with some approachable questions that further explore the relationship
between these techniques and whose answers may assist in future
advances towards a classification scheme.

\end{abstract}

\maketitle

\section{Introduction}

Despite spectacular advances in defining invariants for
simply-connected smooth and symplectic 4-dimensional manifolds and the
discovery of important qualitative features about these manifolds, we
seem to be retreating from any hope to classify simply-connected
smooth or symplectic 4-dimensional manifolds. The subject is rich in
examples that demonstrate a wide variety of disparate phenomena. Yet
it is precisely this richness which, at the time of these lectures,
gives us little hope to even conjecture a classification scheme. In
these notes, adapted from two talks given at the 2004 Clay Institute
Summer School on {\it Floer homology, gauge theory, and low
dimensional topology} at the Alfred R\'enyi Institute, we will quickly
review what we do and do not know about the existence and uniqueness
of smooth and symplectic structures on closed, simply-connected
$4$-manifolds. We will then list the techniques used to date and
capture the key features common to all these techniques. We finish
with some approachable questions that further explore the relationship
between these techniques and whose answers may assist in future
advances towards a classification scheme.

\subsection*{Algebraic Topology} 
The critical algebraic topological information for a closed,
simply-connected, smooth $4$-manifold $X$ is encoded in its Euler
characteristic $e(X)$, its signature $\s(X)$, and its type $t(X)$
(either $0$ if the intersection form of $X$ is even and $1$ if it is
odd). These invariants completely classify the homeomorphism type of
$X$ (\cites{donaldson,freedman}).  We recast these algebraic
topological invariants by defining $\ch(X)=(e(X)+\s(X))/4$, which is
the holomorphic Euler characteristic in the case that $X$ is a complex
surface, and $c(X)=3\s(X)+2e(X)$, which is the self-intersection of
the first Chern class of $X$ in the case that $X$ is complex.

\subsection*{Analysis} 
To date, the critical analytical information for a smooth, closed,
simply-connected $4$-manifold $X$ is encoded in its Seiberg-Witten
invariants \cite{witten}. If $b_2^+(X)>1$, this integer-valued
function $\ssw_X$ is defined on the set of $spin ^c$ structures over
$X$. Corresponding to each $spin ^c$ structure $\fs $ over $X$ is the
bundle of positive spinors $W^+_{\fs}$ over $X$. Set $c(\fs)\in
H_2(X)$ to be the Poincar\'e dual of $c_1(W^+_{\fs})$.  Each $c(\fs)$
is a characteristic element of $H_2(X;\Z)$ (i.e. its Poincar\'e dual
$\hat{c}(\fs)=c_1(W^+_{\fs})$ reduces mod~2 to $w_2(X)$).  The sign of
$\ssw_X$ depends on a homology orientation of $X$, that is, an
orientation of $H^0(X;\R)\otimes\det H_+^2(X;\R)\otimes \det
H^1(X;\R)$. If $\ssw_X(\b)\neq 0$, then $\b$ is called a {\it{basic
class}} of $X$. It is a fundamental fact that the set of basic classes
is finite. Furthermore, if $\b$ is a basic class, then so is $-\b$
with $\ssw_X(-\b)=(-1)^{(\it{e}+\text{sign})(X)/4}\,\ssw_X(\b)$
where $\text{e}(X)$ is the Euler number and $\text{sign}(X)$ is the
signature of $X$. The Seiberg-Witten invariant is an
orientation-preserving diffeomorphism invariant of $X$ (together with
the choice of a homology orientation). We recast the Seiberg-Witten
invariant as an element of the integral group ring $\Z H_2(X)$, where
for each $\a\in H_2(X)$ we let $t_\a$ denote the corresponding element
in $\Z H_2(X)$. Suppose that $\{\pm \b_1,\dots,\pm \b_n\}$ is the set
of nonzero basic classes for $X$.  Then the Seiberg-Witten invariant
of $X$ is the Laurent polynomial
\[\sw_X = \ssw_X(0)+\sum_{j=1}^n \ssw_X(\b_j)\cdot
(t_{\b_j}+(-1)^{(\it{e}+\text{sign})(X)/4}\, t_{\b_j}^{-1}) \in \Z H_2(X).\]

In the case $\ch =1$, the Seiberg-Witten invariant depends on a given orientation of $H^2_+(X;\R)$,  a given metric $g$, and a self-dual 2-form as follows. 
There is a unique $g$-self-dual harmonic 2-form $\o_g\in H^2_+(X;\R)$ with $\o_g^2=1$ and corresponding to the positive orientation. Fix a characteristic homology class $k\in H_2(X;\Z)$.  Given a pair $(A,\psi)$, where
$A$ is a connection in the complex line bundle whose first Chern class is the Poincar\'e dual $\hk=\frac{i}{2\pi}[F_A]$ of $k$ and $\psi$ a section of the bundle $W^+$ of self-dual spinors for the associated $spin^{\,c}$ structure, the perturbed Seiberg-Witten equations are:
\begin{gather*} 
D_A\psi = 0 \\
F_A^+  = q(\psi)+i\eta \notag\label{SWeqn}
\end{gather*}
where $F_A^+$ is the self-dual part of the curvature $F_A$,
$D_A$ is the twisted Dirac operator, $\eta$ is a
self-dual 2-form on $X$, and
$q$ is a quadratic function. Write $\ssw_{X,g,\eta}(k)$ for the
corresponding invariant. As the pair
$(g,\eta)$ varies, $\ssw_{X,g,\eta}(k)$ can change only at those pairs
$(g,\eta)$ for which there are solutions with $\psi=0$. These 
solutions occur for pairs $(g,\eta)$ satisfying $(2\pi\hk+\eta)\cdot\o_g=0$.
This last equation defines a wall in $H^2(X;\R)$. 

The point $\o_g$ determines a component of the double cone consisting of elements of $H^2(X;\R)$ of positive square. We prefer to work with $H_2(X;\R)$. The dual component is determined by the Poincar\'e dual $H$ of $\omega_g$.  (An element $H'\in H_2(X;\R)$ of positive square lies in the same component as $H$ if $H'\cdot H>0$.) If
$(2\pi \hk+\eta)\cdot\o_g\ne 0$ for a generic $\eta$, $\,\ssw_{X,g,\eta}(k)$ is
well-defined, and its value depends only on the sign of $(2\pi \hk+\eta)\cdot\o_g$. Write $\ssw_{X,H}^+(k)$ for $\ssw_{X,g,\eta}(k)$ if 
$(2\pi \hk+\eta)\cdot\o_g>0$ and $\ssw_{X,H}^-(k)$ in the other case.

The invariant $\ssw_{X,H}(k)$ is defined by $\ssw_{X,H}(k) =\ssw_{X,H}^+(k)$ if 
$(2\pi \hk)\cdot\o_g>0$, or dually, if $k\cdot H>0$, and $\ssw_{X,H}(k) =\ssw_{X,H}^-(k)$ if $H\cdot k <0$. The wall-crossing formula \cites{KMgenus,LiLiu} states that if $H', H''$ are elements of positive square in $H_2(X;\R)$ with $H'\cdot H>0$ and $H''\cdot H>0$, then if $k\cdot H' <0$ and $k\cdot H''>0$,
\[ \ssw_{X,H''}(k) - \ssw_{X,H'}(k) = (-1)^{1+\frac12 d(k)}\]
where $d(k)=\frac14(k^2-(3\,\text{sign}+2\,\it{e})(X))$ is the formal dimension of the Seiberg-Witten moduli spaces.

Furthermore, in case $b^-\le 9$, the wall-crossing formula, together with the fact that $\ssw_{X,H}(k)=0$ if $d(k)<0$, implies that $\ssw_{X,H}(k) = \ssw_{X,H'}(k)$ for any $H'$ of positive square in $H_2(X;\R)$ with $H\cdot H'>0$. So in case $b^+_X=1$ and $b^-_X\le 9$, there is a well-defined Seiberg-Witten invariant, $\ssw_X(k)$.

\subsection*{Possible Classification Schemes} From this point forward and unless otherwise stated all manifolds will be closed and simply-connected. In order to avoid trivial constructions we consider {\it irreducible} manifolds, i.e. those that cannot be represented as the connected sum of two manifolds except if one factor is a homotopy $4$-sphere. (We still do not know if there exist smooth homotopy $4$-spheres not diffeomorphic to the standard $4$-sphere $S^{4}$).

So the existence part of a classification scheme for irreducible smooth (symplectic) $4$-manifolds could take the form of determining which $(\ch, c, t) \in \Z \times \Z \times \Z_{2}$ can occur as $(\ch(X),c(X), t(X))$ for some smooth (symplectic) $4$-manifold $X$. This is referred to as {\it the geography problem}. The game plan would be to create techniques to realize all possible lattice points. The uniqueness part of the classification scheme would then be to determine all smooth (symplectic) $4$-manifolds with a fixed $(\ch(X),c(X), t(X))$ and determine invariants that would distinguish them. Again, the game plan would be to create techniques that preserve the homeomorphism type yet change these invariants.  

In the next two sections we will outline what is and is not known about the existence (geography) and uniqueness problems without detailing the techniques. Then we will list the techniques used, determine their interplay, and explore questions that  may yield new insight.  A companion approach, which we will also discuss towards the end of these lectures,  is to start with a particular well-understood class of $4$-manifolds and determine how all other smooth (symplectic) $4$-manifolds can be  constructed from these.

\section {Existence}

Our current understanding of the geography problem is given by Figure 1 where all known simply-connected smooth irreducible $4$-manifolds are plotted as lattice points in the $(\ch, c)$-plane. In particular, all known simply-connected irreducible smooth or symplectic $4$-manifolds have $ 0 \le c < 9\ch $ and every lattice point in that region can be realized by a symplectic (hence smooth) $4$-manifold.

\setlength{\unitlength}{1in}
\begin{picture}(4.5,4.5)
\put(.2,.5){\vector(0,1){3.5}}
\put(.2,.5){\vector(1,0){4}}
\put(.2,.5){\line(1,3){1.15}}
\put(.2,.5){\line(1,2){1.7}}
\put(.8,.5){\line(3,2){3.4}}
\put(.8,.5){\line(3,1){3.4}}
\put(0,3.4){$c$}
\put(4,.35){$\ch$}
\put(.3,.25){Elliptic Surfaces $E(n)$  ($(\ch,c)=(n,0)$)}
\put(3.1,.3){$ c < 0$  ??}
\put(1.4,2.6){surfaces of general type}
\put(1.5,2.35){$2\ch - 6 \le c \le 9\ch$}
\put(1.1,4.1){$c=9\ch$}
\put(.3,3.1){$c>9\ch$ ??}
\put(1.8,4.1){$c=8\ch$}
\put(3.9,2.8){$c=2\ch - 6$}
\put(3.9,1.41){$c=\ch - 3$}
\put(1.15,3.1){$\sigma >0$}
\put(2.25,3.1){$\sigma <0$}
\put(3.2,2.0){symplectic with}
\put(3.1,1.85){one $SW$ basic class}
\put(3.0,1.70){$\ch - 3 \le c \le 2\ch -6$}
\put(3.2,1.50){(cf. \cite{FPS})}
\put(3.1,1.1){symplectic with}
\put(2.77,.95){$(\ch-c-2)$  $SW$ basic classes}
\put(3.1,.80){$0 \le c \le (\ch-3)$}
\put(3.3,.63){(cf. \cite{FScan})}
\multiput(.170,.475)(.2,0){19}{$\bullet$}
\put (1.8, 0){Figure 1}
\end{picture}

An irreducible $4$-manifold need not lie on a lattice point. The issue here is whether $\ch \in \Z$ or  $\ch \in \Z[\frac{1}{2}]$. Note that $\ch(X) \in \Z$ iff $X$ has an almost-complex structure. In addition, the Seiberg-Witten invariants are only defined for manifolds with $\ch \in \Z$. To date our only technique to determine if a manifold is irreducible is to use Seiberg-Witten invariants. Thus, all known results have $\ch \in \Z$.

\begin{prob} Do there exist irreducible smooth 4-manifolds with $\ch \notin \Z$?
\end{prob}

Here the work of Bauer and Furuta \cite{BF} on stable homotopy invariants derived from the Seiberg-Witten equations may be useful. To expose our ignorance, consider two copies of the elliptic surface $E(2)$. Remove the neighborhood of a sphere with self-intersection $-2$ from each and glue together the resulting manifolds along their boundary ${{\mathbf{RP}}^{\,3}}$ using the orientation reversing diffeomorphism of ${{\mathbf{RP}}^{\,3}}$. The result has $\ch \notin \Z$ and it is unknown if it is irreducible.

All complex manifolds with $c=9\ch > 9$ are non-simply-connected, in particular they are ball quotients.  
Thus obvious problems are:

\begin{prob} Do there exist irreducible simply connected smooth or symplectic manifolds with $c=9\ch > 9$?
\end{prob}

\begin{prob} Does there exist an irreducible non-complex smooth  or symplectic manifold $X$  with  $\ch >1$, $c=9\ch$ (with any fundamental group), $\sw_X \ne 0$, and  which is not a ball-quotient?
\end{prob}

\begin{prob} Do there exist irreducible smooth or symplectic manifolds with $c> 9\ch$?
\end{prob}

The work of Taubes \cite{TGW} on the relationship between Seiberg-Witten and Gromov-Witten invariants shows that $c \ge 0$ for an irreducible symplectic $4$-manifold.

\begin{prob} Do there exist simply connected irreducible smooth  manifolds with $c < 0$?
\end{prob}

There appears to be an interesting relationship between the number of Seiberg-Witten basic classes and the pair $(\ch,c)$. In particular, all known smooth $4$-manifolds with $0  \le c \le \ch-3$ have at least $\ch - c -2 $ Seiberg-Witten basic classes \cite{FPS}. So

\begin{prob} Does there exist an irreducible smooth  manifold $X$  with $0  \le c(X) \le \ch(X)-3$ and with fewer than $\ch(X)-c(X)-2$ Seiberg-Witten basic classes? (There is a physics proof that there are no such examples \cite{marino}.)
\end{prob}

Figure 1 contains no information about the geography of spin $4$-manifolds, i.e. manifolds with $t=0$. For a spin $4$-manifold there is the relation $c = 8\ch \mod 16$. Almost every lattice point with $c = 8\ch \mod 16$ and $0 \le c < 9\ch$ can be be realized by an irreducible spin 4-manifold \cite{Park}. Surprisingly not all of the lattice points with $2\ch \le 3(\ch - 5) $ can be realized by complex manifolds \cite{persson}, so these provide several examples of smooth irreducible $4$-manifolds with $2\ch-6 \le c < 9 \ch$ that support no complex structure (cf. \cite{families}). The open issues here are a better understanding of manifolds close to the $c=9\ch$ line, in particular those with $9\ch > c \ge 8.76\ch$ and not on the lines $c=9\ch - k$ with $k \le 121$ (cf. \cite{persson}).

The techniques used in all these constructions are an artful application of the generalized fiber sum  construction (cf. \cite{GS})  and the rational blowdown construction \cite{rat}, which we will  discuss later in this lecture.

\section{Uniqueness}

Here is where we begin to lose control of the classification of smooth $4$-manifolds.  If a topological $4$-manifold admits an irreducible smooth (symplectic) structure that has a smoothly (symplectically) embedded torus with self-intersection zero and with simply-connected complement, then it also admits infinitely many distinct smooth (symplectic) structures and also admits infinitely many distinct smooth structures with no compatible symplectic structure. The basic technique here is the knot-surgery construction of Fintushel-Stern \cite{KL4M}, i.e. remove a neighborhood $T^{2} \times D^{2}=S^{1}\times S^{1}\times D^{2}$ of this torus and replace it with $S^{1}\times S^{3}\setminus K$ where $K$ is a knot in $S^{3}$. In particular, there are no known examples of (simply-connected) smooth or symplectic $4$-manifolds with $\chi >1$ that does not admit such a torus. Hence, there are no known smooth or symplectic $4$-manifolds with $\ch > 1$ that admit finitely many smooth or symplectic structures. Thus,
\begin{prob} Do there exist irreducible smooth (symplectic) 4-manifolds with $\ch >1$ that do not  admit a smoothly (symplectically) embedded torus with self-intersection 0 and simply-connected complement? 
\end{prob}

All of the constructions used for the geography problem  with $\ch > 1$ naturally contain such tori, so the only hope is to find manifolds where these constructions have yet to work, i.e. those  with $8.76 < c \le 9\ch$, that do not contain such tori. 

\begin{prob} Do manifolds with $c = 9\ch$ admit exotic smooth structures?
\end{prob}

The situation for $\ch =1$ is potentially more interesting and may
yield phenomena not shared by manifolds with $\ch >1$.  For example,
the complex projective plane $\CP$ has $c=9\ch = 9$ and is
simply-connected. It is also known that $\CP$ as a smooth manifold has a unique symplectic
structure \cites{taubes, TGW}. Thus, a fundamental question that still
remains is

\begin{prob} Does the complex projective plane $\CP$ admit exotic smooth structures?
\end{prob}

\begin{prob} What is the smallest $m$ for which $\CP\# m\,\CPb$ admits an exotic smooth structure?
\end{prob}

The primary reason for our ignorance here is  that for $c>1$ (i.e. $m < 9$), these manifolds do not contain homologically essential tori with zero self-intersction. Since the rational elliptic surface $E(1)\cong \CP\# 9\,\CPb$ admits tori with self-intersection zero, it has infinitely many distinct smooth structures.  In the late 1980's Dieter Kotschick   \cite{K} proved that the Barlow surface, which was known to be homeomorphic to $\CP\#\, 8\CPb$, is not diffeomorphic to it. In the following years the subject of simply connected smooth $4$-manifolds with $m<8$  languished because of a lack of suitable examples.   However, largely due to a beautiful paper of Jongil Park \cite{Park2}, who found the first examples of exotic simply connected $4$-manifolds with $m=7$, interest was revived. Shortly after this conference ended, Peter Ozsv\'ath and Zolt\'an Szab\'o proved that Park's manifold is minimal \cite{OS} by computing its Seiberg-Witten invariants. Then Andr\'as Stipsicz and Zolt\'an Szab\'o used a technique similar to Park's to construct an exotic manifold with $m=6$ \cite{SS}. The underlying technique in these constructions is an artful use of the rational blowdown construction.

Since $\CP\# m\,\CPb$ for $m<9$ does not contain smoothly embedded tori with self-intersection zero, it has not been known whether it can have an infinite family of smooth structures. Most recently, Fintushel and Stern \cite{FSfake} introduced a new technique which was used to show that for $6 \le m \le 8$,  $\CP\# m\,\CPb$ does indeed have an infinite family of smooth structures, and, in addition,  none of these smooth structures support a compatible symplectic structure. Park, Stipsicz, and Szab\'o \cite{PSS} used this construction to show that $m=5$ also has an infinite family of smooth structures none of which support a compatible symplectic structure (cf.  \cite{FSfake}). The basic technique in these constructions is a prudent blend of the knot surgery and  rational blowdown constructions.

As is pointed out in \cite{SS}, the Seiberg-Witten invariants will never distinguish more than two distinct irreducible symplectic structures on  $\CP\# m\,\CPb$ for $m<9$. Basically, this is due to the fact that if there is more than one pair of basic classes for a $\ch=1$ manifold, then it is not minimal. So herein lies one of our best hopes for finiteness in dimension $4$.

\begin{prob} Does $\CP\# m\,\CPb$ for $m<9$ support more than two irreducible symplectic structures that are not deformation equivalent? 
\end{prob}

\section{The techniques used for the construction of all known  simply-connected smooth and symplectic $4$-manifolds}

The construction of simply-connected smooth or symplectic
$4$-manifolds sometimes takes the form of art rather than
science. This is exposed by the lack of success in proving structural
theorems or uncovering any finite phenomena in dimension 4.  In this
lecture we will list all the constructions used in building the
$4$-manifolds necessary for the results of the first two sections and
try to bring all the unusual phenomena in dimension 4 into a framework
that will allow us to at least understand those surgical operations
that one performs to go from one smooth structure on a given
simply-connected 4-manifold to any other smooth structure. This will
take the form of understanding a variety of cobordisms between
4-manifolds.

Here is the list of constructions used in the first two sections.

\begin{description}
\item[generalized fiber sum] Assume two $4$-manifolds $X_{1}$ and
$X_{2}$ each contain an embedded genus $g$ surface $F_{j} \subset
X_{j}$ with self-intersection $0$. Identify tubular neighborhoods $\nu
F_{j}$ of $F_{j}$ with $F_{j}\times D^{2}$ and fix a diffeomorphism
$f:F_{1} \to F_{2}$. Then the fiber sum $X= X_{1}\#_{f} X_{2}$ of
$(X_{1}, F_{1})$ and $(X_{2},F_{2})$ is defined as $X_{1}\setminus \nu
F_{1} \cup_{\phi}X_{2}\setminus \nu F_{2}$, where $\phi$ is $f \times$
(complex conjugation) on the boundary $\partial(X_{j}\setminus \nu
F_{j})=F_{j}\times S^{1}$.

\item[generalized logarithmic transform] Assume that $X$ contains a
homologically essential torus $T$ with self-intersection zero. Let
$\nu T$ denote a tubular neighborhood of $T$. Deleting the interior of
$\nu T$ and regluing $T^{2} \times D^{2}$ via a diffeomorphism $\phi:
T^{2} \times D^{2} \to \partial (X- {\rm int }\ \nu T)= \partial\nu T$
we obtain a new manifold $X_{\phi}$, the generalized logarithmic
transform of $X$ along $T$.

If $p$ denotes the absolute value of the degree of the map $\pi \circ
\phi: \{pt\}\times S^{1}\to \pi(\partial\nu T)=S^{1}$, then $X_{\phi}$
is called a generalized logarithmic transformation of multiplicity
$p$.

If the complement of $T$ is simply-connected and $t(X)=1$, then
$X_{\phi}$ is homeomorphic to $X$. If the complement of $T$ is
simply-connected and $t(X)=0$, then $X_{\phi}$ is homeomorphic to $X$
if $p$ is odd, otherwise $X_{\phi}$ has the same $c$ and $\ch$ but
with $t(X_{\phi})=1$.
\item[blowup] Form $X \#\CPb$.
\item[rational blowdown ] Let $C_{p}$ be the smooth $4$-manifold
obtained by plumbing $(p-1)$ disk bundles over the $2$-sphere
according to the diagram

\setlength{\unitlength}{.013in}
\begin{picture}(100,60)(-90,-25)
 \put(-12,3){\makebox(200,20)[bl]{$-(p+2)$ \hspace{6pt}
                                  $-2$ \hspace{96pt} $-2$}}
 \put(4,-25){\makebox(200,20)[tl]{$u_{0}$ \hspace{25pt}
                                  $u_{1}$ \hspace{86pt} $u_{p-2}$}}
  \multiput(10,0)(40,0){2}{\line(1,0){40}}
  \multiput(10,0)(40,0){2}{\circle*{3}}
  \multiput(100,0)(5,0){4}{\makebox(0,0){$\cdots$}}
  \put(125,0){\line(1,0){40}}
  \put(165,0){\circle*{3}}
\end{picture}

\noindent Then the classes of the $0$-sections have self-intersections
$u_0^2=-(p+2)$ and $u_i^2=-2$, $i=1,\dots,p-2$. The boundary of $C_p$
is the lens space $L(p^2, 1-p)$ which bounds a rational ball $B_p$
with $\pi_1(B_p)=\Z_p$ and $\pi_1(\bd B_p)\to \pi_1(B_p)$
surjective. If $C_p$ is embedded in a $4$-manifold $X$ then the
rational blowdown manifold $X_{(p)}$ is obtained by replacing $C_p$
with $B_p$, i.e., $X_{(p)} = (X\- C_p) \cup B_p$ (cf. \cite{rat}).
If $X\- C_p$ is simply connected, then so is $X_{(p)}$

\item[knot surgery] Let $X$ be a $4$-manifold which contains a
homologically essential torus $T$ of self-intersection $0$, and let
$K$ be a knot in $S^3$. Let $N(K)$ be a tubular neighborhood of $K$ in
$S^3$, and let $T\x D^2$ be a tubular neighborhood of $T$ in $X$. Then
the knot surgery manifold $X_K$ is defined by
\[ X_K = (X\- (T\x D^2))\cup (S^1\x (S^3\- N(K))\]
The two pieces are glued together in such a way that the homology
class $[{\text{pt}}\x \bd D^2]$ is identified with $[{\text{pt}}\x
\lam]$ where $\lam$ is the class of a longitude of $K$.  If the
complement of $T$ in $X$ is simply-connected, then $X_{K}$ is
homeomorphic to $X$.
\end{description}

The Seiberg-Witten invariants are sensitive to all of these operations.

\begin{description}
\item[generalized logarithmic transform] If $T$ is contained in a node neighborhood, then 

\[\sw_{X_{\phi}}=\sw_{X}\cdot(s^{-(p-1)}+s^{-(p-3)}+
\cdots+s^{(p-1)})\]
where $s=\exp(T/p)$, $p$ the order of the generalized logarithmic transform (cf. \cite{rat}).
\item[blowup] The relationship between the Seiberg-Witten invariants
of $X$ and its blowup $X \# \CPb$ is referred to as the blowup formula
and was given in Witten's original article \cite{witten}
(cf. \cite{Turkey}). In particular, if $e$ is the homology class of
the exceptional curve and $\{B_1, \dots, B_n\}$ are the basic classes
of $X$, then the basic classes of $X \# \CPb$ are $\{B_1 \pm E, \dots,
B_n \pm E\}$ and $\ssw_{X \# \CPb}(B_j \pm E)=\ssw_X(B_j)$.
\item[rational blowdown] The Seiberg-Witten invariants of $X$ and
$X_{(p)}$ can be compared as follows. The homology of $X_{(p)}$ can be
identified with the orthogonal complement of the classes $u_i$,
$i=0,\dots,p-2$ in $H_2(X;\Z)$, and then each characteristic element
$k\in H_2(X_{(p)};\Z)$ has a lift $\widetilde{k} \in H_2(X;\Z)$ which
is characteristic and for which the dimensions of moduli spaces agree,
$d_{X_{(p)}}(k)=d_X(\widetilde{k})$. It is proved in \cite{rat} that
if $b^+_X>1$ then $\ssw_{X_{(p)}}(k)=\ssw_X(\widetilde{k})$. In case
$b^+_X=1$, if $H\in H_2^+(X;\R)$ is orthogonal to all the $u_i$ then
it also can be viewed as an element of $H_2^+(X_{(p)};\R)$, and
$\ssw_{X_{(p)},H}(k)=\ssw_{X,H}(\widetilde{k})$.

\item[knot surgery]  If, for example, $T$ is contained in a  a node neighborhood and $\ch(X)>1$ then the Seiberg-Witten invariant of  the knot
surgery manifold $X_K$ is given by
\[ \sw_{X_K}=\sw_{X}\cdot\DD_K(t) \]
where $\DD_K(t)$ is the symmetrized Alexander polynomial of $K$ and
$t=\exp(2[T])$. When $\ch=1$, the Seiberg-Witten invariants of $X_{K}$
are still completely determined by those of $X$ and the Alexander
polynomial $\DD_K(t)$ \cite{KL4M}.
\end{description}

\noindent Here $T$ contained in a node neighborhood means that an essential loop on $\partial \nu T$ bounds a disk in the complement with relative self-intersection $-1$. We sometimes refer to this disk as a vanishing cycle.

In many circumstances, there are formulas for determining the Seiberg-Witten invariants of a fiber sum in terms of the Seiberg-Witten invariants of $X_{1}$ and $X_{2}$ and how the basic classes intersect the surfaces $F_{1}$ and $F_{2}$.

\subsection*{Interaction of the operations}

While knot surgery appears to be a new operation, the constructions in
\cite{KL4M} point out that the knot surgery construction is actually a
series of $\pm 1$ generalized logarithmic transformations on {\it
null-homologous tori}. To see this, note that any knot can be
unknotted via a sequence of crossing changes, which in turn can be
realized as a sequence of $\pm 1$ surgeries on unknotted curves
$\{c_1, \dots, c_n\}$ that link the knot algebraically zero times and
geometrically twice. When crossed with $S^{1}$ this translates to the
fact that $X$ can be obtained from $X_{K}$ via a sequence of $\pm 1$
generalized logarithmic transformations on the {\it null-homologous
tori} $\{S^1 \times c_1, \dots, S^1\times c_n\}$ in $X_K$. So the
hidden mechanism behind the knot surgery construction is generalized
logarithmic transformations on null-homologous tori.  The calculation
of the Seiberg-Witten invariants is then reduced to understanding how
the Seiberg-Witten invariants change under a generalized logarithmic
transformation on a null-homologous torus. This important formula is
due to Morgan, Mrowka, and Szab\'o \cite{MMS} (see also \cite{T}). For
this formula fix simple loops $\a$, $\b$, $\d$ on $\bd N(T)$ whose
homology classes generate $H_1(\bd N(T))$. If $\o = p\a +q\b + r\d$
write $X_T(p,q,r)$ instead of $X_T(\o)$.  Given a class $k\in H_2(X)$:

\begin{multline}\label{surgery formula} \sum_i\ssw_{X_T(p,q,r)}(k_{(p,q,r)}+2i[T])= 
p\sum_i\ssw_{X_T(1,0,0)}(k_{(1,0,0)}+2i[T]) +\\+q\sum_i\ssw_{X_T(0,1,0)}(k_{(0,1,0)}+2i[T])
+r\sum_i\ssw_{X_T(0,0,1)}(k_{(0,0,1)}+2i[T])
\end{multline}
In this formula, $T$ denotes the torus which is the core $T^2\x {0}\C T^2\x D^2$ in each specific manifold $X(a,b,c)$ in the formula, and $k_{(a,b,c)}\in H_2(X_T(a,b,c))$ is any class which agrees with the restriction of $k$ in $H_2(X\- T\x D^2,\bd)$ in the diagram:

\[ \begin{array}{ccc}
H_2(X_T(a,b,c)) &\longrightarrow & H_2(X_T(a,b,c), T\x D^2)\\
&&\Big\downarrow \cong\\
&&H_2(X\- T\x D^2,\bd)\\
&&\Big\uparrow \cong\\
H_2(X)&\longrightarrow & H_2(X,T\x D^2)
\end{array}\]

Let $\pi(a,b,c): H_2(X_T(a,b,c))\to H_2(X\- T\x D^2,\bd)$ be the composition of maps in the above diagram, and $\pi(a,b,c)_*$ the induced map of integral group rings. Since we are often interested in invariants of the pair $(X,T)$, it is sometimes useful to work with \[\swb_{(X_T(a,b,c),T)}=\pi(a,b,c)_*(\sw_{X_T(a,b,c)})\in \Z H_2(X\- T\x D^2,\bd).\]
The indeterminacy due to the sum in \eqref{surgery formula} is caused by multiples of $[T]$; so
passing to $\swb$ removes this indeterminacy, and the Morgan-Mrowka-Szab\'o formula becomes
\begin{equation}\label{surgery formula 2}
 \swb_{(X_T(p,q,r),T)} = 
p\swb_{(X_T(1,0,0),T)} +q\swb_{(X_T(0,1,0),T)}
+r\swb_{(X_T(0,0,1),T)}.
\end{equation}

So if we expand the notion of generalized logarithmic transformation to include both homologically essential and null-homologous tori, then we can eliminate the knot surgery construction from our list of essential surgery operations.  Thus our list is of essential operations is reduced to

\begin{itemize}
\item generalized fiber sum
\item generalized logarithmic transformations on a torus with trivial normal bundle
\item blowup
\item rational blowdown
\end{itemize}

There are further relationships between these operations. In
\cite{rat} it is shown that if $T$ is contained in a node
neighborhood, then a generalized logarithmic transformation can be
obtained via a sequence of blowups and rational blowdowns. (This
together with work of Margaret Symington \cite{sym} shows that
logarithmic transformations ($p \ne 0$) on a symplectic torus results
in a symplectic manifold. We do not know of any other proof that a
generalized logarithmic transformations on a symplectic torus in a
node neighborhood results in a symplectic manifold.) However, it is
not clear that a rational blowdown is always the result of blowups and
logarithmic transforms.

Rational blowdown changes the topology of the manifold $X$; while
$\chi$ remains the same, $c$ is decreased by $p-3$. So, an obvious
problem would be

\begin{prob}\label{log}
Are any two homeomorphic simply-connected smooth 4-manifolds related
via a sequence of generalized logarithmic transforms on tori?
\end{prob}

As already pointed out, there are two cases.

\begin{enumerate}
\item $T$ is essential in homology
\item $T$ is null-homologous
\end{enumerate}

This leads to:

\begin{prob}\label{tori}
Can a  generalized logarithmic transform on a homologicaly essential torus be obtained via a sequence of generalized logarithmic transforms on null-homologous tori?
\end{prob}

For the rest of the lecture we will discuss these last two problems.

\section{Cobordisms between $4$-manifolds}.

Let $X\sb{1}$ and $X\sb{2}$ be two homeomorphic simply-connected
smooth 4-manifolds. Early results of C.T.C. Wall show that there is an
$h$--cobordism $W^5$ between $X_1$ and $X_2$ obtained from $X\sb{1}
\times I$ by attaching $n$ $2-$handles and $n$ $3-$handles. A long
standing problem that still remains open is:

\begin{prob}  Can $W^5$ can be chosen so that $n=1$. 
\end{prob}

Let's explore the consequences if we can assume $n=1$. We can then
describe the $h-$cobordism $W^5$ as follows. First, let $W\sb{1}$ be
the cobordism from $X\sb{1}$ to $X\sb{1} \# S\sp{2}\times S\sp{2}$
given by attaching the $2-$handle to $X_1$. To complete $W^5$ we then
would add the $3-$handle. Dually, this is equivalent to attaching a
$2-$handle to $X_2$. So let $W\sb{2}$ be the cobordism from $X\sb{2}$
to $X\sb{2} \# S\sp{2}\times S\sp{2}$ given by attaching this
$2-$handle to $X_2$.  Then $W^5=W\sb{1} \cup\sb{f} (-W\sb{2})$ for a
suitable diffeomorphism $f\colon X\sb{1} \#S\sp{2}\times S\sp{2}
\rightarrow X\sb{2} \# S\sp{2}\times S\sp{2}$. Let $A$ be any of the
standard spheres in $S\sp{2}\times S\sp{2}$. Then the complexity of
the $h$-cobordism can be measured by the type $k$, which is half the
minimum of the number of intersection points between $A$ and $f(A)$
(as $A \cdot f(A)=0$ there are $k$ positive intersection points and
$k$ negative intersection points). This complexity has been studied in
~\cite{MS}. A key observation is that if $k=1$, then a neighborhood of
$A \cup f(A)$ is diffeomorphic to an embedding of twin spheres in
$S^4$ and that its boundary is the three-torus $T^3$. A further
observation is that $X_2$ is then obtained from $X_1$ by removing a
neighborhood of a null-homologous torus $T$ embedded in $X_1$ (with
trivial normal bundle) and sewing it back in differently. Thus when
$k=1$, $X_2$ is obtained from $X_1$ by a generalized generalized
logarithmic transform on a {\it null-homologous} torus.

Thus, the answers to Problems~\ref{log} and \ref{tori} are clearly
related to the complexity $k$ of $h$-cobordisms. We expect that the
answer to Problem~\ref{tori} is NO and that ordinary generalized
logarithmic transforms on homologically essential tori will provide
examples of homeomorphic $X_1$ and $X_2$ that require $h$--cobordisms with
arbitrarily large complexity.

Independent of this, an important next step is to study complexity
$k>2$ $h$-cobordisms. Here, new surgical techniques are suggested. In
particular, the neighborhood of $A \cup f(A)$ above is diffeomorphic
to the neighborhood $N^{\prime}$ of two $2-$spheres embedded in $S^4$
with $2n$ points of intersection. Let $N$ be obtained from
$N^{\prime}$ with one of the $2-$spheres surgered out. Then it can be
shown that $X_1$ is obtained from $X_2$ by removing an embedding of
$N^{\prime}$ and regluing along a diffeomorphism of its boundary. This
could lead to a useful generalization of logarithmic transforms along
null-homologous tori. It would then be important to compute its effect
on the Seiberg-Witten invariants, and reinterpret generalized
logarithmic transforms from this point of view.

\subsection*{Round handlebody cobordisms}

Suppose that $X_1$ and $X_2$ are two manifolds with the same $c$ and $\ch$. It follows from early work of Asimov ~\cite{asimov} that there is a round handlebody cobordism $W$ between $X_1$ and $X_2$.  Thus $X_1$ can be obtained from $X_2$ by attaching a sequence of round $1-$handles and round $2-$handles. A round handle is just $S^1$ times a handle in one lower dimension. So for us, a round $r-$handle is a copy of $S^1\times (D^{r}\times D^{4-r})$ attached along $S^1\times (S^{r-1}\times D^{4-r})$ (see ~\cite{asimov} for definitions).

\begin{prob}\label{rd}
Can $W$ be chosen so that there are no round $1-$handles?
\end{prob}

For a moment, suppose that the answer to Problem~\ref{rd} is Yes. Then
$W$ would consist of only round $2-$handles. It then follows that
$X_2$ would be obtained from $X_1$ via a sequence of generalized
logarithmic transforms on tori. Thus the answer to Problem~\ref{rd} is
tightly related to Problem~\ref{log}.

Note that if $X_1$ and $X_2$ are round handlebody cobordant, then the only invariant preventing them from being homeomorphic is whether $t(X_{1})=t(X_{2})$. So suppose $t(X_{1})=0$ and $t(X_2)=1$. If the answer to Problem~\ref{rd} were yes, then one could change the second Stiefel-Whitney class via a sequence of generalized logarithmic transforms on tori. By necessity these tori cannot be null-homologous. So understanding new surgical operations that will change $t$  without changing $c$, $\ch$, and preserving the Seiberg-Witten invariants should provide new insights.

\begin{prob}\label{rd1} Suppose two simply-connected smooth 4-manifolds have the  same $c$, $\ch$, number of Seiberg-Witten basic classes, and different $t$. Determine surgical operations that will transform one to the other.
\end{prob}

There are explicit examples of this phenomena amongst complex surfaces, e.g. two Horikawa surfaces with the same $c$ and $\ch$, but different $t$. 

\section*{Modifying symplectic $4$-manifolds}

To finish up this lecture, we point out that all known constructions of (simply-connected) non-symplectic $4$-manifolds can be obtained from symplectic $4$-manifolds by performing logarithmic transforms on null-homologous Lagrangian tori with non-vanishing framing defect (cf. \cite{Lagr}). Let's look at a specific example of this phenomena. In particular, let's consider $E(n)_K$.

The elliptic surface $E(n)$ is the double branched cover of $\SS$
with branch set equal to
four disjoint copies of $S^2\x \{{\rm{pt}}\}$ together with
$2n$ disjoint copies of $\{{\rm{pt}}\}\x S^2$. The resultant branched
cover has $8n$ singular
points (corresponding to the double points in the branch set), whose
neighborhoods are cones
on
${{\mathbf{RP}}^{\,3}}$. These are desingularized in the usual way,
replacing their
neighborhoods with cotangent bundles of $S^2$. The result is $E(n)$.
The horizontal and
vertical fibrations of
$\SS$ pull back to give fibrations of
$E(n)$ over $\PO$. A generic fiber of the vertical fibration is the
double cover of $S^2$,
branched over $4$ points --- a torus. This describes an elliptic fibration of
$E(n)$. The generic fiber of the horizontal fibration is the double cover of
$S^2$, branched over $2n$ points, and this gives a genus $n-1$ fibration on
$E(n)$. This genus $n-1$ fibration has four singular fibers which are
the preimages of the
four $S^2\x \{{\rm{pt}}\}$'s in the branch set  together with the
spheres of self-intersection
$-2$ arising from desingularization.  The generic fiber
$T$ of the elliptic fibration meets a generic fiber
$\Sig_{n-1}$ of the horizontal fibration in two points,
$\Sigma_{n-1}\cdot T=2$.

Now let $K$ be a fibered knot of genus $g$, and fix a generic elliptic
fiber $T_0$ of
$E(n)$. Then in the knot surgery manifold
\[ E(n)_K = (E(n)\- (T_0\x D^2)) \cup (S^1\x (S^3\- N(K)) ,\] each
normal $2$-disk to $T_0$
is replaced by a fiber of the fibration of $S^3\- N(K)$ over
$S^1$. Since
$T_0$ intersects each generic horizontal fiber twice, we obtain a
`horizontal' fibration
\[ h: E(n)_K\to\PO\]  of genus $2g+n-1$.

This fibration also has four singular fibers arising from the four
copies of $S^2\x
\{{\rm{pt}}\}$ in the branch set of the double cover of $\SS$. Each
of these gets blown up at
$2n$ points in $E(n)$, and the singular fibers each  consist of a genus
$g$ surface $\Sig_g$ of self-intersection $-n$ and multiplicity $2$
with $2n$ disjoint
$2$-spheres of self-intersection $-2$, each meeting $\Sig_g$
transversely in one point. The
monodromy around each singular fiber is (conjugate to) the diffeomorphism of
$\Sig_{2g+n-1}$ which is the deck transformation $\eta$ of the double cover of
$\Sig_g$, branched over $2n$ points. Another way to describe $\eta$
is to take the
hyperelliptic involution $\o$ of $\Sig_{n-1}$ and to connect sum copies of
$\Sig_g$ at the two points of a nontrivial orbit of $\o$. Then $\o$
extends to the involution
$\eta$ of $\Sig_{2g+n-1}$.

The fibration which we have described is not Lefschetz since the
singularities are not simple
nodes. However, it can be perturbed locally to be Lefschetz. 

So in summary, if $K$ is a fibered knot whose fiber has genus
$g$, then $E(n)_K$  admits a locally holomorphic fibration (over
$\PO$) of genus $2g+n-1$ which has exactly four singular fibers.
Furthermore, this fibration can be deformed locally to be Lefschetz.

There is another way
to view these constructions. Consider the branched double cover of $\Sig_g\x S^2$ whose branch set
consists of two disjoint copies of $\Sig_g\x \{{\rm{pt}}\}$ and $2n$
disjoint copies of
$\{{\rm{pt}}\}\x S^2$. After desingularizing as above, one obtains a complex surface denoted
$M(n,g)$. Once again, this  manifold carries a pair of
fibrations. There is a genus
$2g+n-1$ fibration over $S^2$  and an $S^2$ fibration over $\Sig_g$.

Consider first the $S^2$ fibration. This has $2n$ singular fibers,
each of which consists of
a smooth 2-sphere $E_i$, $i=1,\dots, 2n$, of  self-intersection $-1$
and multiplicity $2$,
together with a pair of disjoint spheres of self-intersection
$-2$, each intersecting $E_i$ once transversely. If we blow down
$E_i$ we obtain again an $S^2$ fibration over $\Sig_g$, but the $i$th
singular fiber  now
consists of a pair of 2-spheres of self-intersection $-1$ meeting
once, transversely. Blowing
down one of these gives another $S^2$ fibration over $\Sig_g$, with
one less singular fiber.
Thus blowing down $M(n,g)$ $4n$ times results in  a manifold which is
an $S^2$ bundle over
$\Sig_g$. This shows that (if $n>0$) $M(n,g)$ is diffeomorphic to
$(S^2\x \Sig_g)\# 4n\,\CPb$.

The genus $2g+n-1$ fibration on $M(n,g)$ has 2 singular fibers. As
above, these fibers consist
of a genus $g$ surface $\Sig_g$ of self-intersection $-n$ and
multiplicity $2$ with $2n$
disjoint $2$-spheres of self-intersection $-2$, each meeting $\Sig_g$
transversely in one
point. The monodromy of the  fibration around each of these fibers is
the deck transformation
of the double branched cover of
$\Sig_g$. This is just the map $\eta$ described above.

Let $\vp$ be a diffeomorphism of $\Sig_g\- D^2$ which is the identity
on the boundary. For
instance, $\vp$ could be the monodromy of a fibered knot of genus
$g$. There is an induced diffeomorphism $\Phi$ of \,
       $\Sig_{2g+n-1}=\Sig_g\#\Sig_{n-1}\#\Sig_g$ which is given by $\vp$
on the first
$\Sig_g$ summand and by the identity on the other summands. Consider
the twisted fiber sum
\[ M(n,g)\#_{\Phi}M(n,g) = \{M(n,g)\- (D^2\x \Sig_{2g+n-1})\}
\cup_{{\text{id}}\x\Phi} \{M(n,g)\- (D^2\x \Sig_{2g+n-1})\} \] where
fibered neighborhoods of
generic fibers $\Sig_{2g+n-1}$ have been removed from the two copies
of $M(n,g)$, and they
have been glued by the diffeomorphism
${{\text{id}}\x\Phi}$ of $S^1\x \Sig_{2g+n-1}$.

In the case that  $\vp$ is the monodromy of a fibered knot $K$, it can be shown that
$M(n,g)\#_{\Phi}M(n,g)$ is the manifold $E(n)_K$ with the genus
$2g+n-1$ fibration described above. To see this, we view $S^2$ as the
base of the horizontal
fibration. Then it suffices to check that the total monodromy map
$\pi_1(S^2\- 4\, {\text{points}})\to {\text{Diff}}(\Sig_{2g+n-1})$ is
the same for  each. It
is not difficult to see that if we write the generators of
$\pi_1(S^2\- 4\, {\text{points}})$ as $\a$, $\b$, $\g$ with $\a$ and
$\b$ representing loops
around the singular points of, say, the image of  the first copy of
$M(n,g)$ and basepoint in
this image, and $\g$ a loop around a  singular point in the image of
the second $M(n,g)$ then
the monodromy map $\mu$ satisfies
$\mu(\a)=\eta$, $\mu(\b)=\eta$ and $\mu(\g)$ is
$\vp\oplus\o\oplus\vp^{-1}$, expressed as a
diffeomorphism of $\Sig_g\#\Sig_{n-1}\#\Sig_g$. That  this is also
the monodromy of $E(n)_K$
follows directly from its construction.

Now let $E(n)_g$ denote $E(n)$ fiber summed with $T^2\times\Sig_g$
along an elliptic fiber.  The penultimate observation is that
$E(n)_K$, viewed as $M(n,g)\#_{\Phi}M(n,g)$, is then the result of a
sequence of generalized logarithmic transforms on null-homologous
Lagrangian tori in $E(n)_g$.  The effect of these surgeries is to
change the monodromy of the genus $n+2g-1$ Lefschetz fibration (over
$\PO$) on $E(n)_g$. This is accomplished by doing a $1/n$, with
respect to the natural Lagrangian framing, generalized logarithmic
transform on these Lagrangian tori (cf. \cites{families, Lagr}). The
final observation is that if the Lagrangian framing of these tori
differs from the null-homologous framing (cf. ~\cite{Lagr}), then a
$1/n$ log transformations on $T$ with respect to the null-homologous
framing can be shown, by computing Seiberg-Witten invariants, to
result in non-symplectic 4-manifolds. Careful choices of these tori
and framings will result in manifolds homotopy equivalent to
$M(n,g)\#_{\Phi}M(n,g)$ (cf. \cite{families}).

\end{document}